\documentclass[aop,MSNbibl,dvips]{arximspdf}
\usepackage{mathbh}

\doi{10.1214/12-AOP753} 
\volume{41}
\issue{5}
\pubyear{2013}
\firstpage{3518}
\lastpage{3541}

\makeatletter
\newcommand{\rrvert}{\vert}

\newtheorem{theorem}{Theorem}
\newtheorem{lem}[theorem]{Lemma}
\newtheorem{prop}[theorem]{Proposition}

\newproclaim{rmk}{Remark}

\newcommand{\Q}{\mathbb{Q}}
\renewcommand{\P}{\mathbb{P}}
\newcommand{\E}{\mathbb{E}}
\newcommand{\R}{\mathbb{R}}
\newcommand{\F}{\mathcal{F}}
\newcommand{\ind}{\mathbh{1}}

\makeatother

\begin{document}
\begin{frontmatter}

\title{A simple path to asymptotics for the frontier of a~branching
Brownian motion}
\runtitle{Asymptotics for the frontier of BBM}

\begin{aug}
\author[A]{\fnms{Matthew I.} \snm{Roberts}\corref{}\ead[label=e1]{mattiroberts@gmail.com}}
\runauthor{M. I. Roberts}
\affiliation{Weierstrass Institute for Applied Analysis and Stochastics}
\address[A]{Weierstrass Institute\\
\quad for Applied Analysis\\
\quad and Stochastics\\
Mohrenstrasse 39\\
10117 Berlin\\
Germany\\
\printead{e1}} 
\end{aug}

\received{\smonth{10} \syear{2011}}
\revised{\smonth{2} \syear{2012}}

%
\begin{abstract}
We give short proofs of two classical results about the position of the
extremal particle in a branching Brownian motion, one concerning the
median position and another the almost sure behaviour.
\end{abstract}

%

\end{frontmatter}

\section{Introduction and main results}\label{intro}

Kolmogorov et al.~\cite{kolmogorovetaletudekpp} proved that the
extremal particle in a standard branching Brownian motion sits near
$\sqrt2 t$ at time~$t$. Higher order corrections to this result were
given by Bramson~\cite{bramsonmaximaldisplacementBBM}, and then
almost sure fluctuations were proved by Hu and Shi
\cite{hushiminimalposcritmgconvBRW}. These two remarkable papers,
more than thirty years apart, provide results which reflect an
extremely deep understanding of the underlying branching structure.
This article grew out of a desire to know whether shorter or simpler
proofs of these results exist.

We consider a branching Brownian motion (BBM) beginning with one
particle at 0, which moves like a standard Brownian motion until an
independent exponentially distributed time with parameter 1. At this
time it dies and is replaced (in its current position) by two new
particles, which---relative to their birth time and position---behave
like independent copies of their parent, moving like Brownian motions
and branching at rate 1 into two copies of themselves. Let $N(t)$ be
the set of all particles alive at time~$t$, and if $v\in N(t)$, then
let $X_v(t)$ be the position of $v$ at time $t$. If $v\in N(t)$ and
$s<t$, then let $X_v(s)$ be the position of the unique ancestor of $v$
that was alive at time $s$. Define $M_t = \max_{v\in N(t)} X_v(t)$.

\subsection{Bramson's result on the distribution of $M_t$}

Define
\[
u(t,x) = \P( M_t \leq x ).
\]
This function $u$ satisfies the
Fisher--Kolmogorov--Petrovski--Piscounov equation
\[
u_t = \tfrac{1}{2}u_{xx} + u^2 - u
\]
(with Heaviside initial condition), which has been studied for many
years both analytically and probabilistically; see, for example,
Kolmogorov et al.~\cite{kolmogorovetaletudekpp},
Fisher~\cite{fisheradvanceadvantageousgenes}, Skorohod
\cite{skorohodbranchingdiffusionprocesses}, McKean
\cite{mckeanapplicationbbmkpp}, Bramson
\cite{bramsonmaximaldisplacementBBM,bramsonconvergenceKoleqntravwaves}, Neveu
\cite{neveumultiplicativemartingalesspatialbps}, Uchiyama
\cite{uchiyamanonlineardiffusionequationslargetime}, Aronson and
Weinberger~\cite{aronsonweinbergernonlineardiffusionsingenetics},
Karpelevich et al.
\cite{karpelevichetalbranchingdiffusionskpptravwaves}, Harris~\cite{harrissctravwavesfkppprobarguments},
Kyprianou
\cite{kyprianoualternativesimonprobanalysisfkpp}, Harris et al.
\cite{harrisetalfkpponesidedtravelingwaves}. In particular (see
\cite{kolmogorovetaletudekpp}) $u$ converges to a \textit{travelling
wave}: that is, there exist functions $m$ of $t$ and $w$ of $x$ such
that $1-w$ is a probability distribution function, and
\[
u\bigl(t, m(t)+x\bigr)\to w(x)
\]
uniformly in $x$ as $t\to\infty$.

We note that $m$ and $w$ are not uniquely determined by this
definition; since we shall be concerned with the detailed behaviour of
$m$, to be precise we set $m(t):= \sup\{x\in\R\dvtx\P(M_t\leq x) \leq
1/2\}$. We offer a proof of the following result which is shorter and
simpler than the original proof by Bramson
\cite{bramsonmaximaldisplacementBBM}:

%
\begin{theorem}[(Bramson~\cite{bramsonmaximaldisplacementBBM})]\label
{bramsonthm}
The median $m(t)$ satisfies
\[
m(t) = \sqrt{2}t - \tfrac{3}{2\sqrt2}\log t + O(1) \qquad\mbox{as }
t\to\infty.
\]
\end{theorem}

As Bramson notes in~\cite{bramsonmaximaldisplacementBBM}, ``an
immediate frontal assault using moment estimates, but ignoring the
branching structure of the process, will fail.'' That is, for $y\geq0$
define $\beta= \sqrt2 - \frac{3}{2\sqrt2}\frac{\log t}{t} + \frac
{y}{t}$ and let $G(t)$ be the set of particles near $\beta t$ at time
$t$. If some particle has large position at time $s<t$, then many
particles are likely to have large position at time $t$, and therefore
the moments of $\#G(t)$ are misleading.

To get around this problem we consider a subset $H(t)$ of $G(t)$. A
lower bound for $m(t)$ will follow if we can show that the first two
moments of $\#H(t)$ are well behaved. Our approach differs from
Bramson's only in that our set $H(t)$ is simpler than his, being the
set of particles that stay below the straight line $\beta s + 1$ for
all $s\leq t$ and end near $\beta t$. This drastically reduces the
difficulty of the calculations required for bounding the moments and is
one reason why our proof is much shorter than the original.

For the upper bound we are forced to return to a more complicated set
$\Gamma(t)$ which is the set of particles that stay below a carefully
chosen curve \mbox{$f(s)+y+1$}, $s\leq t$, and end near $\beta t$. Calculation
of $\E[\#\Gamma(t)]$ is more difficult than that of $\E[\#H(t)]$,
but the two quantities turn out to be of roughly the same size. The key
observation now is that if a particle reaches $f(s)+y$ for some $s<t$,
then it has done the hard work and is likely to have descendants near
$\beta t$, even if we insist that they stay below $f(r)+y+1$ for all
$r\in[s,t]$. Thus the probability that some particle reached $f(s)+y$
for some $s<t$ cannot be much larger than $\E[\#\Gamma(t)]$.

\subsection{Hu and Shi's result on the paths of $M_t$}

Having established Bramson's result on the centring term $m(t)$, we
move on to the almost sure behaviour of $M_t$. We prove the following
result, which is the analogue of a result for quite general branching
random walks given by Hu and Shi~\cite{hushiminimalposcritmgconvBRW}.

%
\begin{theorem}\label{hushithm}
The maximum $M_t$ satisfies
%
%
\begin{equation}
\label{convliminf}
\liminf_{t\to\infty} \frac{M_t - \sqrt{2}t}{\log t} = -\frac
{3}{2\sqrt2} \qquad\mbox{almost surely}
\end{equation}
and
\begin{equation}
\label{convlimsup}
\limsup_{t\to\infty} \frac{M_t - \sqrt{2}t}{\log t} = -\frac
{1}{2\sqrt2} \qquad\mbox{almost surely}.
\end{equation}
\end{theorem}

Thus, although by Theorem~\ref{bramsonthm} the extremal particle
looks like $m(t)$ for most times $t$, occasionally a particle will
travel much further from the origin. Technically the theorem as stated
here is a new result as Hu and Shi considered only discrete-time
branching random walks, but it would not take too much effort to derive
it from their work. We proceed instead by directly applying the
estimates developed in the proof of Theorem~\ref{bramsonthm}. Only
the lower bound in (\ref{convlimsup}) requires a significant amount
of extra work, and for that we take an approach similar to that of Hu
and Shi in~\cite{hushiminimalposcritmgconvBRW}. They noticed
that although the probability that a particle has position much bigger
than $m(t)$ at a fixed time $t$ is very small, the probability that
there exists a time $s$ between (say) $n$ and $2n$ such that a particle
has position much bigger than $m(s)$ at time $s$ is actually quite
large. Here we again simplify the calculations by considering the
number of particles staying below a straight line rather than a curve,
much as in our lower bound for Theorem~\ref{bramsonthm}.

\subsection{Extensions and other models}

We note that although we consider only the simplest possible BBM, with
binary branching at fixed rate 1, our methods can be applied to rather
more general models. There is, however, one important necessary
condition for the proof of our lower bound, that \textit{the mean and
variance of the number of particles born at a branching event must be
finite}. This is simply due to the fact that we employ a second moment method.

Addario-Berry and Reed~\cite{addarioberryreedminimabrw} (in their
Theorem 3) proved an analogue of Bramson's result (our Theorem \ref
{bramsonthm}) for a wide class of branching random walks. We
conjecture that the ideas presented in this article could also be used
to give a new proof of the Addario-Berry and Reed result, relaxing the
conditions on bounded family sizes and independence amongst families.
However, this task would require substantial extra technical work. The
estimates on Bessel processes used to estimate numbers of particles
staying below straight lines can be replaced by small deviations
probabilities for random walks conditioned to stay positive (see
\cite{vatutinwachtellocalprobs}); but calculating the expected number of
particles staying below a curved line, our Lemma~\ref{gammalem},
becomes much more difficult; see the footnote on page 756 of
\cite{hushiminimalposcritmgconvBRW}. Finally, one must make sure that
particles do not jump too far beyond this curved line, which can be
done with conceptually standard but technically delicate first moment estimates.

In a sense, Bramson~\cite{bramsonconvergenceKoleqntravwaves} improved
the $O(1)$ error in Theorem~\ref{bramsonthm}, showing that under his
definition one could choose $m(t)$ such that the corresponding error
was $o(1)$. A related result for branching random walks was recently
given by A\"id\'ekon~\cite{aidekonconvergencelawbrw}, showing
convergence to a specified law for the recentred extremal particle.
This extra detail requires new ideas and is beyond the scope of our
methods.

\subsection{Notation}

We will often use positive constants $c_1,c_2,\ldots$ that are
independent of all other parameters. We shall reset the subscripts at
the end of each proof, so the $c_1$ appearing in the proof of Lemma
\ref{besslem} is not necessarily the same constant as the $c_1$
appearing in Lemma~\ref{besslem2}, for example. On the other hand,
$C_1,C_2,\ldots$ will be positive constants that are fixed throughout
the article.

\section{Bessel-3 processes}\label{besssec}
We begin by recalling some very basic properties of Bessel-3 processes.
If $W_t$, $t\geq0$, is a Brownian motion in $\R^3$ started from
$(x,0,0)$, then its modulus $|W_t|$, $t\geq0$, is called a Bessel-3
process started from $x$. For aesthetic purposes in this article we
shall simply write ``Bessel process'' when we mean ``Bessel-3
process.'' Suppose that $B_t$ is a Brownian motion in $\R$ started from
$B_0 = x$ under a probability measure $P_x$; then $X_t:=
x^{-1}B_t\ind_{\{B_s > 0\ \forall s\leq t\}}$ is a nonnegative
unit-mean martingale under $P_x$. We may change measure by $X_t$,
defining a new probability measure $\hat P_x$ via
\[
\frac{d\hat P_x}{dP_x}\bigg|_{\F_t}:= X_t
\]
(where $\F_t$ is the natural filtration of the Brownian motion $B_t$)
and then $B_t$, $t\geq0$, is a Bessel process under $\hat P_x$. The
density of a Bessel process satisfies
\[
\hat P_x(B_t \in dz) = \frac{z}{x\sqrt{2\pi t}}
\bigl(e^{-(z-x)^2/2t} - e^{-(z+x)^2/2t} \bigr) \,dz.
\]
This and much more about Bessel processes can be found in many
textbooks, for example, Revuz and Yor~\cite{revuzyorctsmartingalesbm}.

%
\begin{lem}\label{bessbound}
Let $\gamma= 2^{1/2}/\pi^{1/2}$. For any $t>0$ and $x,z\geq0$,
\[
\frac{\gamma z^2}{t^{3/2}} e^{-z^2/2t - x^2/2t} \leq\frac{z}{x\sqrt
{2\pi t}} \bigl(e^{-(z-x)^2/2t}
- e^{-(z+x)^2/2t} \bigr) \leq\frac
{\gamma z^2}{t^{3/2}}.
\]
\end{lem}

\begin{pf}
The lower bound is trivial since
\[
e^{xz/t} - e^{-xz/t} = 2\sinh(xz/t) \geq2xz/t.
\]
For the upper bound, note that
\begin{eqnarray*}
\frac{d}{dz} \bigl(e^{-{(z-x)^2}/({2t})} - e^{-{(z+x)^2}/({2t})} \bigr)
&=&
\frac{x}{t} \bigl(e^{-{(z-x)^2}/({2t})} + e^{-{(z+x)^2}/({2t})}
\bigr)\\
&&{} +
\frac{z}{t} \bigl(e^{-{(z+x)^2}/({2t})} - e^{-{(z-x)^2}/({2t})} \bigr) \\
&\leq&
\frac{2x}{t}
\end{eqnarray*}
so $e^{-(z-x)^2/2t} - e^{-(z+x)^2/2t} \leq2xz/t$.
\end{pf}

The two lemmas that follow do much of the dirty work of Theorem \ref
{bramsonthm} and Proposition~\ref{existio} (which is the most
difficult part of Theorem~\ref{hushithm}) by calculating the
expectation of two functionals of two dependent Bessel-3 processes.
These calculations will not be motivated until later in the article,
but we include them here as they are facts about Bessel processes that
do not contribute a great deal to the main ideas of the proofs. We
start with Lemma~\ref{besslem}, which will be used in proving the
lower bound for Theorem~\ref{bramsonthm}.

Suppose that under $\hat P$ we have two processes $Y^1_t$ and $Y^2_t$,
$t\geq0$, and a time $\tau\in[0,\infty)$ such that:
\begin{itemize}
\item($Y^1_t$, $t\geq0$) is a Bessel process started from $1$;
\item$\tau$ is exponentially distributed with parameter $2$, and is
independent of $(Y^1_t$, \mbox{$t\geq0)$};
\item$Y^2_t = Y^1_t$ for all $t\leq\tau$;
\item conditioned on $\tau$ and $(Y^1_t, t\leq\tau)$,
$(Y^2_{t+\tau}, t\geq0)$ is a Bessel processes started from
$Y^1_\tau$ that is independent of $(Y^1_t, t>\tau)$.
\end{itemize}
It is clear from this description that $(\tau, Y^1_\tau, Y^1_t,
Y^2_t)$ has a well-behaved joint density. Note that we continue to use
$\hat P$ for this setup, as well as for the single Bessel process
$(B_t, t\geq0)$ seen above.

%
\begin{lem}\label{besslem}
Let
\begin{eqnarray*}
\beta&=& \sqrt{2} - \frac{3}{2\sqrt2}\frac{\log t}{t} + \frac{y}{t},
\\
A_1 &=& \bigl\{1\leq Y^1_t \leq2\bigr\} \quad\mbox
{and}\quad A_2 = \bigl\{1\leq Y^2_t\leq2\bigr\}.
\end{eqnarray*}
There exists a constant $C_1$ such that for all $y\geq0$ and large $t$,
\[
\hat P \bigl[ Y^1_\tau e^{2\tau-({3\tau\log t})/({2t})- \beta
Y^1_\tau}
\ind_{A_1 \cap A_2\cap\{\tau\leq t\}} \bigr] \leq C_1 t^{-3}.
\]
\end{lem}

\begin{pf}
We use the density of $\tau$ to rewrite
\begin{eqnarray*}
&&\hat P \bigl[ Y^1_\tau e^{2\tau-({3\tau\log t})/({2t})- \beta
Y^1_\tau}
\ind_{A_1 \cap A_2\cap\{\tau\leq t\}} \bigr]\\
&&\qquad = 2\int_0^t \hat P
\bigl[Y^1_s e^{-({3s\log t})/({2t}) - \beta Y^1_s}\ind_{A_1 \cap A_2}|
\tau=s
\bigr]\,ds.
\end{eqnarray*}
The idea then is that the probability that a Bessel process is near the
origin at time $t$ is approximately $t^{-3/2}$. If $s$ is small, then
we have two (almost) independent Bessel processes which must both be
near the origin at time $t$, giving $t^{-3}$. If $s$ is large, then we
effectively have only one Bessel process, giving $t^{-3/2}$, but the
$\exp(\frac{3\log t}{2t}s)$ gives us an extra $t^{-3/2}$. When $s$ is
neither large nor small, the above effects combine so that things turn
out nicely. In each case we apply the upper bound from Lemma~\ref{bessbound}.

We first check the small $s$ case:
\begin{eqnarray*}
&&\int_0^1\hat P \bigl[ Y^1_s
e^{-({3s\log t})/({2t}) - \beta
Y^1_s}\ind_{A_1 \cap A_2}| \tau=s \bigr]\,ds
\\
&&\qquad\leq\int_0^1 \hat P(A_1 \cap
A_2 | \tau=s) \,ds
\\
&&\qquad\leq\int_0^1 \hat P \biggl[ \int
_0^\infty\hat P \bigl(Y^1_t,
Y^2_t\in[1,2] | \tau=s, Y^1_s=x
\bigr)\hat P\bigl(Y^1_s\in dx\bigr)\Big\rrvert\tau=s \biggr]
\,ds
\\
&&\qquad\leq\int_0^1 \hat P \biggl[\int
_0^\infty\biggl(\int_1^2
\frac{2z^2}{\sqrt{2\pi}(t-s)^{3/2}}\,dz \biggr)^2 \hat P\bigl(Y^1_s
\in dx\bigr)\Big|\tau=s \biggr] \,ds
\\
&&\qquad\leq c_1 t^{-3},
\end{eqnarray*}
where the third inequality uses Lemma~\ref{bessbound}.
For the large $s$ case,
\[
\int_{t-1}^t \hat P \bigl[ Y^1_s
e^{-({3s\log t})/({2t}) -
\beta Y^1_s}\ind_{A_1 \cap A_2}| \tau=s \bigr]\,ds \leq c_2
t^{-3/2} \hat P(A_1) \leq c_3 t^{-3},
\]
where we have used the fact that, since $\beta\geq1$, $xe^{-\beta
x}\leq1$. (We will use the fact that $\beta\geq1$ throughout the
article without further mention.) Finally the main case, for $s\in[1,t-1]$,
\begin{eqnarray*}
&&\int_1^{t-1}\hat P \bigl[Y^1_s
e^{-({3s\log t})/({2t}) -
\beta Y^1_s}\ind_{A_1 \cap A_2}| \tau=s \bigr]\,ds
\\
&&\qquad\leq\int_1^{t-1} \int_0^\infty
\frac{z^3}{s^{3/2}} e^{- \beta z -
({3s\log t})/({2t})} \biggl(\int_1^2
\frac{2x^2}{\sqrt{2\pi
(t-s)^3}}\,dx \biggr)^2 \,dz \,ds
\\
&&\qquad\leq c_4\int_1^{t-1}
\frac{e^{- ({3s\log t})/({2t})}}{s^{3/2}(t-s)^3} \int_0^\infty z^3
e^{-z} \,dz \,ds
\\
&&\qquad\leq c_5 \int_1^{t-1}
\frac{e^{- ({3s\log t})/({2t})}}{s^{3/2}(t-s)^3} \,ds,
\end{eqnarray*}
where for the first inequality we applied Lemma~\ref{bessbound}. It
is a simple task to bound the last integral above by $t^{-3}$ times a constant
\begin{eqnarray*}
\int_1^{t-1} \frac{e^{- ({3s\log t})/({2t})}}{s^{3/2}(t-s)^3} \,ds
&\leq&
\frac{c_6}{t^3}\int_1^{2t/3} \frac{1}{s^{3/2}}
\,ds + \frac
{c_7}{t^3}\int_{2t/3}^{t-\sqrt t}
e^{-\log t} \,ds\\
&&{} + \frac
{c_8}{t^3}\int_{t-\sqrt t}^{t-1}
\frac{e^{({3\log t})/({2\sqrt
t})}}{(t-s)^3}\,ds \\
&\leq& c_9 t^{-3},
\end{eqnarray*}
which completes the proof.
\end{pf}

Our next lemma is very similar; it estimates a slightly different
functional, which will appear in Proposition~\ref{existio} (the most
difficult part of Theorem~\ref{hushithm}).

%
\begin{lem}\label{besslem2}
Let $\beta_t = \sqrt2 - \frac{1}{2\sqrt2}\frac{\log t}{t}$ and
$a_{s,t} = \frac{1}{2\sqrt2}\log s - \frac{1}{2\sqrt2}\frac{\log
t}{t}s$. If $e\leq s\leq t\leq2s$, then
\begin{eqnarray*}
&&
\hat P \bigl[Y^1_\tau e^{2\tau- ({\tau\log t})/({2t})-
\beta_t Y^1_\tau}
\ind_{\{a_{s,t} + 1 \leq Y^1_s \leq a_{s,t} + 2\}} \ind_{\{1 \leq Y^2_t
\leq2\}} \ind_{\{\tau\leq s\}} \bigr]
\\
&&\qquad\leq C_2 e^{-({s\log t})/({2t})} \biggl(\frac{1}{t^{5/2}} +
\frac{1}{t^{3/2}(t-s+1)^{3/2}} \biggr)
\end{eqnarray*}
for some constant $C_2$ not depending on $s$ or $t$.
\end{lem}

\begin{pf}
Just as in the proof of Lemma~\ref{besslem}, we use the density of
$\tau$ to rewrite
\begin{eqnarray*}
&&
\hat P \bigl[Y^1_\tau e^{2\tau- ({\tau\log t})/({2t})-
\beta_t Y^1_\tau}
\ind_{\{a_{s,t} + 1 \leq Y^1_s \leq a_{s,t} + 2\}} \ind_{\{1 \leq Y^2_t
\leq2\}} \ind_{\{\tau\leq s\}} \bigr]
\\
&&\qquad
=2\int_0^s e^{-({r\log t})/({2t})} \hat P \bigl[
Y^1_r e^{-\beta_t Y^1_r}\ind_{\{a_{s,t} + 1 \leq Y^1_s \leq a_{s,t} +
2\}}
\ind_{\{1 \leq Y^2_t \leq2\}} | \tau=r \bigr] \,dr
\end{eqnarray*}
and then approximate the integral. Essentially the $e^{-\beta_t
Y^1_r}$ term means our initial Bessel process must be near the origin
at time $r$; then two independent Bessel processes started from time
$r$ must be near the origin at times $s$ and $t$, respectively. If $r\in
[1,s-1]$, then integrating out over $Y^1_r$, applying Lemma~\ref
{bessbound} three times and using the fact that $\int_0^\infty z^3
e^{-\beta_t z}\,dz < \infty$,
\begin{eqnarray*}
&&
\hat P \bigl[Y^1_r e^{-\beta_t Y^1_r}
\ind_{\{a_{s,t} + 1 \leq
Y^1_s \leq a_{s,t} + 2\}} \ind_{\{1 \leq Y^2_t \leq2\}} | \tau=r \bigr
] \\
&&\qquad\leq c_1
\int_0^\infty z e^{-\beta_t z}
\frac
{z^2}{r^{3/2}} \cdot\frac{1}{(s-r)^{3/2}} \cdot\frac
{1}{(t-r)^{3/2}} \,dz
\\
&&\qquad\leq c_2 r^{-3/2}(s-r)^{-3/2}(t-r)^{-3/2}.
\end{eqnarray*}
For $r\leq1$ we are effectively asking two independent Bessel
processes to be near the origin at times $s$ and $t$, giving
$s^{-3/2}t^{-3/2}$, and for $r\geq s-1$ we have just one Bessel process
which must be near the origin\vadjust{\goodbreak} at times $s$ and $t$, giving
$s^{-3/2}(t-s+1)^{-3/2}$. These two simple calculations follow as in
Lemma~\ref{besslem}. Thus
\begin{eqnarray*}
&&
\int_0^s e^{-({r\log t})/({2t})} \hat P \bigl[
Y^1_r e^{-\beta_t Y^1_r} \ind_{\{a_{s,t} + 1 \leq Y^1_s \leq a_{s,t}
+ 2\}}
\ind_{\{1 \leq Y^2_t \leq2\}} | \tau=r \bigr] \,dr
\\
&&\qquad
\leq\frac{c_3}{s^{3/2} t^{3/2}} + c_4 \int_1^{s-1}
\frac{e^{-({r\log t})/({2t})}}{r^{3/2}(s-r)^{3/2}(t-r)^{3/2}} \,dr + \frac
{c_5 e^{-({s\log t})/({2t})}}{s^{3/2} (t-s+1)^{3/2}}.
\end{eqnarray*}
Since $s$ and $t$ are of the same order, and $\log s \geq\frac{\log
t}{t}s$ provided $s, t \geq e$,
it remains to estimate the integral in the last line above. We proceed
again just as in Lem\-ma~\ref{besslem}. First the small $r$ case,
\begin{eqnarray*}
\int_1^{s/2}
\frac{e^{-({r\log t})/({2t})}}{r^{3/2}(s-r)^{3/2}(t-r)^{3/2}} \,dr &\leq&
c_6\int_1^{s/2} \frac
{1}{r^{3/2}s^{3/2}t^{3/2}}
\,dr \\
&\leq&\frac{c_7}{s^{3/2}t^{3/2}},
\end{eqnarray*}
the large $r$ case,
\[
\int_{s-s/t^{1/4}}^{s-1} \frac{e^{-({r\log t})/({2t})}}{r^{3/2}(s-r)^{3/2}(t-r)^{3/2}} \,dr \leq
\frac{c_8 e^{-({s\log t})/({2t})}}{s^{3/2}(t-s+1)^{3/2}}
\]
and finally the intermediate $r$ case,
\begin{eqnarray*}
\int_{s/2}^{s - s/t^{1/4}} \frac{e^{-({r\log t})/({2t})}}{r^{3/2}(s-r)^{3/2}(t-r)^{3/2}} \,dr &\leq&
c_9 \frac
{t^{3/4}}{s^{9/2}} \int_{s/2}^{s-s/t^{1/4}}
e^{-({r\log t})/({2t})} \,dr
\\
&\leq& c_{10} \frac{t^{7/4}}{s^{9/2}} e^{-
({s\log t})/({4t})} \\
&\leq&
\frac{c_{11}}{t^{5/2}} e^{-({s\log t})/({2t})},
\end{eqnarray*}
where we have again used $\log s \geq\frac{\log t}{t}s$ and $s\leq
t\leq2s$.
\end{pf}

\section{The many-to-one and many-to-two lemmas}\label{manytofewsec}
We mentioned in the \hyperref[intro]{Introduction} that we will attempt
to count the
number of particles remaining below certain lines and ending near
$\beta t$. To do this we will need to calculate the first two moments
of the number of such particles. In this section we state results for
doing so in the form that will be most useful to us. These are standard
first and second moment bounds for branching processes combined with
one- and two-particle changes of measure.

\subsection{The many-to-one lemma}\label{manytoonesec}
The many-to-one lemma is a simple and well-known tool in branching
processes. It essentially says that the expected number of particles
with a certain property equals the expected number of particles times
the probability that one particle has that property. To be more
precise, let\vadjust{\goodbreak} $g_t(v)$ be a measurable functional of $t$ and the path of
a particle $v$ up to time $t$; so, for example, we might take
\[
g_t(v) = \ind_{\{X_v(t)\geq x\}} \quad\mbox{or}\quad g_t(v) =
t^2 e^{\int_0^t X_v(s) \,ds}.
\]
Then the standard many-to-one lemma says
\[
\E\biggl[\sum_{v\in N(t)} g_t(v) \biggr] =
e^t E\bigl[g_t(\xi)\bigr],
\]
where $\xi_t$, $t\geq0$, is just a standard Brownian motion under $P$.

Now, sometimes it will be easiest to calculate $E[g_t(\xi)]$ by using
a change of measure. Fixing $\alpha>0$ and $f\dvtx[0,\infty)\to\R$ such
that $f\in C^2$, and defining
\[
\zeta(t) = \frac{1}{\alpha}\bigl(\alpha+ f(t) - \xi_t
\bigr)e^{\int_0^t
f'(s) \,d\xi_s - (\int_0^t f'(s)^2 \,ds)/2}\ind_{\{\xi_s<\alpha
+ f(s)\ \forall s\leq t\}},
\]
the following lemma is a result of Girsanov's theorem and the knowledge
of Bessel processes at the start of Section~\ref{besssec}. It will be
useful for counting the number of particles near $\beta t$ that have
remained below $\alpha+ f(s)$ for all $s\leq t$. For a proof see
Theorem 8.5 of~\cite{hardyharrisspineapproachapplications}.

%
\begin{lem}[(Many-to-one lemma)]\label{manytoone}
\[
\E\biggl[\sum_{v\in N(t)} g_t(v)
\ind_{\{X_v(s)<\alpha+f(s)\
\forall s\leq t\}} \biggr] = e^t \Q\biggl[\frac{1}{\zeta(t)}
g_t(\xi) \biggr],
\]
where under $\Q$, $\alpha+ f(t) - \xi_t$, $t\geq0$, is a Bessel process.
\end{lem}

\subsection{The many-to-two lemma}
We also use a many-to-two lemma, which---just as the many-to-one lemma
reduces calculating first moments to consideration of just one
particle---will reduce calculating second moments to functionals of
two, necessarily dependent, particles. This is a natural idea, and
Bramson uses a basic many-to-two lemma in
\cite{bramsonmaximaldisplacementBBM}. Again we will combine this idea
with a change of measure. [Note, however, that while we used a general
$C^2$ function $f$ in our many-to-one lemma, we will need only $f(s) =
\beta s$ here.] We do not prove Lemma~\ref{manytotwo}---as Bramson
says, ``a rigorous verification is quite messy''---and refer to Lemma 3
of~\cite{harrisrobertsmultiplespines} which gives a quite general
formulation.

Suppose that under $\Q$, as well as the process $\xi_t$ seen in
Section~\ref{manytoonesec}, we have two processes $\xi^1_t$ and
$\xi^2_t$, $t\geq0$, and a time $T\in[0,\infty)$ such that:
\begin{itemize}
\item$(1 + \beta t - \xi^1_t, t\geq0)$ is a Bessel processes
started from $1$;
\item$T$ is exponentially distributed with parameter $2$, and is
independent of $(\xi^1_t$, \mbox{$t\geq0)$};\vspace*{1pt}
\item$\xi^2_t = \xi^1_t$ for all $t\leq T$;
\item conditioned on $T$ and $(\xi^1_t, t\leq T)$, $(\beta
(T+s) -
\xi^2_{T+s}, s\geq0)$ is a Bessel processes started from
$\beta T
- \xi^1_T$ that is independent of $(\xi^1_t, t>T)$.\vadjust{\goodbreak}
\end{itemize}
Define
\[
\zeta^i(t) = \bigl(1 + \beta t - \xi^i_t
\bigr)e^{\beta\xi^i_t - \beta^2
t/2}\ind_{\{\xi^i_s<1 + \beta s \ \forall s\leq t\}}
\]
for $i=1,2$ and $t\geq0$.

%
\begin{lem}[(Many-to-two lemma)]\label{manytotwo}
Let $g_t(\cdot)$ and $h_t(\cdot)$ be measurable functionals of $t$
and the path of a particle up to time $t$, as in Section \ref
{manytoonesec}. Then
\begin{eqnarray*}
&&\E\biggl[\sum_{u,v\in N(t)} g_t(u)
h_t(v)\ind_{\{X_u(s)< 1+\beta
s\ \forall s\leq t, X_v(s)< 1 + \beta s\ \forall s\leq t\}} \biggr]
\\
&&\qquad= \Q\biggl[e^{2t+T\wedge t}\frac{\zeta^1(T\wedge t)}{\zeta^1(t)\zeta
^2(t)}g_t\bigl(
\xi^1\bigr)h_t\bigl(\xi^2\bigr) \biggr]
\\
&&\qquad= e^{3t}\Q\biggl[\frac{1}{\zeta^1(t)}\ind_{\{T>t\}}
g_t\bigl(\xi^1\bigr)h_t\bigl(
\xi^1\bigr) \biggr]\\
&&\qquad\quad{} + e^{2t}\Q\biggl[e^T
\frac{\zeta^1(T)}{\zeta^1(t)\zeta^2(t)}\ind_{\{T\leq t\}} g_t\bigl(\xi^1
\bigr) h_t\bigl(\xi^2\bigr) \biggr].
\end{eqnarray*}
\end{lem}

The dependence between the two Bessel processes reflects the dependence
structure of the BBM: any pair of particles $(u,v)$ in the BBM are
dependent up until their most recent common ancestor. The first term on
the right-hand side above takes account of the possibility that the
Bessel processes have not yet split (which corresponds to the event
that $u$ and $v$ are in fact the same particle) and otherwise the
second term incorporates the split time $T$ of the two Bessel processes
(which corresponds to the last time at which the most recent common
ancestor of $u$ and $v$ was alive).

\section{\texorpdfstring{Proof of Theorem \protect\ref{bramsonthm}}{Proof of Theorem 1}}

\subsection{\texorpdfstring{The lower bound for Theorem \protect\ref{bramsonthm}}{The lower bound for Theorem 1}}

Fix $t>0$ and set (as in Section~\ref{besssec})
\[
\beta= \sqrt{2} - \frac{3}{2\sqrt2}\frac{\log t}{t} + \frac{y}{t}.
\]
Now define
\[
H(y,t) = \# \bigl\{u\in N(t)\dvtx X_u(s) \leq\beta s + 1 \ \forall s
\leq t, \beta t -1 \leq X_u(t) \leq\beta t \bigr\}.
\]
We shall show that the first two moments of $H(y,t)$ give an accurate
picture of the probability that there is a particle near $\beta t$ at
time $t$. We write $g(y,t) \asymp h(y,t)$ if $c_1 g \leq h \leq c_2 g$
for some strictly positive constants $c_1$ and $c_2$ not depending on
$t$ or $y$.

%
\begin{lem}\label{firstmoment}
For $t\geq1$ and $y\in[0,\sqrt t]$,
\[
\E\bigl[H(y,t)\bigr] \asymp e^{-\sqrt{2}y}.\vadjust{\goodbreak}
\]
\end{lem}

\begin{pf}
We apply the many-to-one lemma with $f(t)=\beta t$ and $\alpha=1$.
\begin{eqnarray*}
\E\bigl[H(y,t)\bigr] &=& e^t \Q\biggl[\frac{1}{\zeta(t)}
\ind_{\{\beta t -
1\leq\xi_t \leq\beta t\}} \biggr]\\
&=& e^t\Q\biggl[\frac{e^{-\beta
\xi_t + \beta^2 t/2}}{\beta t + 1 - \xi_t}
\ind_{\{\beta t - 1\leq
\xi_t \leq\beta t\}} \biggr]
\\
&\asymp& e^{t-\beta^2 t/2} \Q(\beta t - 1\leq\xi_t \leq\beta t)
\\
&\asymp& t^{3/2}e^{-\sqrt{2}y} \Q(1 \leq\beta t + 1 -
\xi_t \leq2).
\end{eqnarray*}
Now, $\beta t + 1 - \xi_t$ is a Bessel process started from $1$ under
$\Q$, so by Lemma~\ref{bessbound},
\[
\Q(1 \leq\beta t + 1 - \xi_t \leq2) \asymp\int_1^2
\frac
{z^2}{t^{3/2}} \,dz \asymp t^{-3/2}.
\]
\upqed
\end{pf}

We now use the second moment of $H(y,t)$ to get a lower bound for $m(t)$.

%
\begin{prop}\label{bramsonlower}
There exists a constant $C_3>0$ such that for $t\geq1$ and $y\in
[0,\sqrt t]$,
\[
\P\bigl(\exists u\in N(t)\dvtx X_u(t) \geq\sqrt{2}t -
\tfrac{3}{2\sqrt
2}\log t + y\bigr) \geq C_3 e^{-\sqrt{2} y}.
\]
\end{prop}

\begin{pf}
By reducing $C_3$ if necessary, it suffices to establish the claim for
all large $t$. For $i=1,2$ let $A'_i = \{\beta t - 1\leq\xi^i_t
\leq\beta t \}$.
By the many-to-two lemma,
\begin{eqnarray*}
\E\bigl[H(y,t)^2\bigr] &=& e^{3t}\Q\biggl[
\frac{\ind_{\{T>t\}}}{\zeta^1(t)}\ind_{A'_1} \biggr] + e^{2t}\Q\biggl[
\frac{e^T\zeta^1(T)}{\zeta^1(t)\zeta^2(t)}\ind_{A'_1\cap A'_2\cap\{
T\leq t\}
} \biggr]
\\
&=& e^t\Q\biggl[\frac{1}{\zeta^1(t)}\ind_{A'_1} \biggr]\\
&&{} +
e^{2t}\Q\biggl[\frac{e^T(\beta T + 1 -\xi^1_T)e^{\beta\xi^1_T - \beta^2
T/2}\ind_{A'_1\cap A'_2\cap\{T\leq t\}}}{(\beta t + 1 - \xi^1_t)(\beta
t + 1 -\xi^2_t) e^{\beta\xi^1_t + \beta\xi^2_t - \beta
^2 t}} \biggr]
\\
&\leq&\E\bigl[H(y,t)\bigr]\\
&&{} + e^{2t - \beta^2 t + 2\beta} \Q\bigl
[e^T\bigl(\beta
T + 1 - \xi^1_T\bigr) e^{\beta\xi^1_T - \beta^2 T/2}
\ind_{A'_1\cap
A'_2\cap\{T\leq t\}} \bigr]
\\
&\leq&\E\bigl[H(y,t)\bigr] \\
&&{} + c_1 t^3 e^{-\sqrt{2}y} \Q
\bigl[\bigl(\beta T + 1 - \xi^1_T\bigr)\\
&&\qquad\hspace*{48.6pt}{}\times e^{2T-({3T\log
t})/({2t}) - \beta(\beta T + 1 -
\xi^1_T)}
\ind_{A'_1\cap A'_2\cap\{T\leq t\}} \bigr],
\end{eqnarray*}
where for the second equality we used that $T$ is an exponential random
variable of parameter 2 independent of the path of $\xi^1$, and for
the final inequality we used that if $y\in[0,\sqrt t]$, then
\[
\beta^2 T = 2T - 3\frac{\log t}{t} T + \frac{2\sqrt2 y}{t}T + O(1).
\]

Under $\Q$, $(\beta s + 1 - \xi^1_s, s\geq0)$ and $(\beta s + 1 -
\xi^2_s, s\geq0)$ are Bessel processes starting from $1$ that are
equal up to $T$ and independent (given $T$ and $\xi^1_T$) after $T$.
Thus, taking notation from Lemma~\ref{besslem}, we have
\[
\E\bigl[H(y,t)^2\bigr] \leq\E\bigl[H(y,t)\bigr] + c_1
t^3 e^{-\sqrt{2}y}\hat P \bigl[Y^1_\tau
e^{2\tau-({3\tau\log t})/({2t})- \beta Y^1_\tau}\ind_{A_1 \cap A_2\cap\{
\tau\leq t\}} \bigr].
\]
Lemma~\ref{besslem} tells us that the $\hat P$-expectation is at most
a constant times $t^{-3}$, so for large $t$ and $y\geq0$,
\[
\E\bigl[H(y,t)^2\bigr] \leq c_2 \E\bigl[H(y,t)\bigr]
\]
for some constant $c_2$ not depending on $y$ or $t$. Using Lemma \ref
{firstmoment} we deduce that
\[
\P\bigl(H(y,t) \neq0\bigr) \geq\frac{\E[H(y,t)]^2}{\E[H(y,t)^2]} \geq c_3
e^{-\sqrt2 y}.
\]
\upqed
\end{pf}

\subsection{\texorpdfstring{The upper bound for Theorem \protect\ref{bramsonthm}}{The upper bound for Theorem 1}}

We use a first moment method for an object similar to $H(y,t)$ together
with an estimate of the probability that a particle ever crosses a
carefully chosen line. Again fix $t$, and define
\[
l(s) = \cases{ \frac{3}{2\sqrt2}\log(s+1), &\quad if $0\leq s\leq t/2$,
\vspace*{2pt}\cr
\frac{3}{2\sqrt2}\log(t-s+1), &\quad if $t/2\leq s \leq t$.}
\]
Unfortunately $l$ is not differentiable at $t/2$, so we now choose a
twice continuously differentiable function $L\dvtx[0,t]\to\R$ such that:
\begin{itemize}
\item$L(s) = l(s)$ for all $s\notin[t/2-1,t/2+1]$;
\item$L(s) = L(t-s)$ for all $s\in[0,t]$;
\item$L''(s)\in[-10/t,0]$ for all $s\in[t/2-1,t/2+1]$.
\end{itemize}
Let $f(s) = \beta s + L(s)$ for $s\in[0,t]$, and define
\[
\Gamma= \#\bigl\{u\in N(t)\dvtx X_u(s)<f(s)+y+1 \ \forall s\leq t,
\beta t - 1 \leq X_u(t) \leq\beta t + y\bigr\}.
\]

%
\begin{lem}\label{gammalem}
There exists $C_4$ such that for all $t\geq1$ and $y\in[0,\sqrt t]$,
\[
\E[\Gamma] \leq C_4 (y+2)^4 e^{-\sqrt2 y}.%
\]
\end{lem}

\begin{pf}
By the many-to-one lemma with $\alpha=y+1$, we have
\[
\E[\Gamma] = e^t\Q\biggl[\frac{y+1}{y+1 + f(t) - \xi_t} e^{-\int
_0^t f'(s) \,d\xi_s + (\int_0^t f'(s)^2 \,ds)/2}
\ind_{\{\beta t
- 1 \leq\xi_t\leq\beta t + y\}} \biggr],
\]
where under $\Q$ the process $y+1 + f(s) - \xi_s$, $s\geq0$, is a
Bessel process. Using the fact that
\[
f'(t)\xi_t = \int_0^t
f'(s)\,d\xi_s + \int_0^t
f''(s)\xi_s \,ds,
\]
which follows from integration by parts, we obtain
\begin{eqnarray*}
\E[\Gamma] &\leq& (y+1)e^t\Q\bigl[e^{-f'(t)\xi_t + \int_0^t
f''(s)\xi_s \,ds + (\int_0^t f'(s)^2 \,ds)/2}
\ind_{\{\xi_t\geq
\beta t-1\}} \bigr]
\\[0.5pt]
& \leq &(y+1)e^t \hat P_{y+1} \biggl[\exp\biggl(-f'(t)\beta t + \int_0^t
f''(s)f(s) \,ds\\[0.5pt]
&&\hspace*{90.7pt}{} + (y+1)\int_0^t f''(s) \,ds\\[0.5pt]
&&\hspace*{90.7pt}{} - \int_0^t f''(s)B_s \,ds +
\frac{1}{2}\int_0^t f'(s)^2 \,ds\biggr)
\ind_{\{B_t \leq y+2\}} \biggr]
\\[0.5pt]
& = &
(y+1)e^{t- \beta^2 t/2 - (\int_0^t L'(s)^2
\,ds)/2} \hat P_{y+1} \bigl[e^{\int_0^t L''(s)(y+1-B_s) \,ds}
\ind_{\{B_t\leq
y+2\}} \bigr],
\end{eqnarray*}
where $(B_s, s\geq0)$ is a Bessel process under $\hat P$. Note that
$t-\frac{1}{2}\beta^2 t = \frac{3}{2}\log t - \sqrt2y + O(1)$, so
\[
\E[\Gamma] \leq c_1(y+1) t^{3/2}e^{-\sqrt2 y} \hat
P_{y+1} \bigl[e^{\int_0^t L''(s)(y+1-B_s) \,ds}\ind_{\{B_t\leq y+2\}}
\bigr].
\]
Now, let
\[
\kappa(s) = \cases{ (s+1)^{2/3}, &\quad if $s\leq t/2$,
\cr
(t-s+1)^{2/3}, &\quad if $s>t/2$;}
\]
then $-\int_0^t L''(s)\kappa(s) \,ds\uparrow\kappa$ for some $\kappa
\in(0,\infty)$. We know that on the event $\{B_t\leq y+2\}$,
$B_s-(y+1)$ will stay well below the curve $\kappa(s)$ with
exceedingly high probability, so the $\hat P_{y+1}$-expectation above
should look like a constant times $(y+2)^3 t^{-3/2}$. The following
calculations verify this fact. We split the event that $B_s-(y+1)$ goes
above $\kappa(s)$ into four possibilities. Either there is a sharp
increase over a small time interval, or $B_s-(y+1)$ is large at some
time of the form $j/t$ for $j\in\mathbb{N}$; in the latter case,
either $(y+1)t^{4/3}\leq j \leq t - (y+1)t^{4/3}$, which is so unlikely
that we can forget about insisting that $B_t\leq y+2$, or $j$ is close
to $0$ or $t^2$, and we may apply the Markov property at time $j/t$.
Indeed, letting $\tilde B_s = (B_s-y-1)/\kappa(s)$,
\begin{eqnarray*}
&&\hat P_{y+1} \bigl[e^{\int_0^t L''(s)(y+1-B_s)\,ds}\ind_{\{B_t\leq
y+2\}} \bigr]
\\[0.5pt]
&&\qquad\leq e^{\kappa}\hat P_{y+1}(B_t\leq y+2)\\[0.5pt]
&&\qquad\quad{} + \sum
_{k=1}^\infty e^{(k+1)\kappa}\hat
P_{y+1} \Bigl(\sup_{s\in[0,t]}\tilde B_s\in[k,k+1
), B_t\leq y+2 \Bigr)
\\
&&\qquad\leq e^{\kappa}(y+2)^3 t^{-3/2}\\
&&\qquad\quad{} + \sum
_{k=1}^\infty e^{(k+1)\kappa
}\sum
_{j=0}^{\lceil t^2\rceil} \hat P_{y+1} \biggl(
\sup_{s\in
[{j}/{t},({j+1})/{t} ]} \tilde B_s \geq(\tilde B_{j/t}\vee
\tilde B_{(j+1)/t}) + \frac{k}{2},\\
&&\qquad\quad\hspace*{247pt} B_t\leq y+2 \biggr)
\\
&&\qquad\quad{} + \sum_{k=1}^\infty e^{(k+1)\kappa}\sum
_{j=1}^{\lceil(y+1)t^{4/3}\rceil}\hat P_{y+1} (\tilde
B_{j/t} \geq k/2, B_t\leq y+2 )
\\
&&\qquad\quad{} + \sum_{k=1}^\infty e^{(k+1)\kappa}\sum
_{j=\lceil
(y+1)t^{4/3}\rceil+1}^{\lceil t^2 - (y+1)t^{4/3}\rceil} \hat P_{y+1}
(\tilde
B_{j/t}\geq k/2 )
\\
&&\qquad\quad{} + \sum_{k=1}^\infty e^{(k+1)\kappa}\sum
_{j=\lceil
t^2 - (y+1)t^{4/3}\rceil+1}^{\lceil t^2\rceil-1}\hat P_{y+1} (\tilde
B_{j/t} \geq k/2, B_t\leq y+2 ).
\end{eqnarray*}
The first double sum is bounded above by
\begin{eqnarray*}
&&
\sum_{k=1}^\infty e^{(k+1)\kappa} \sum
_{j=0}^{\lceil t^2\rceil} \frac{y+2}{y+1}
\P_{y+1} \Bigl(\sup_{s\in[j/t,(j+1)/t]}\tilde B_s \geq(\tilde
B_{j/t}\vee\tilde B_{(j+1)/t}) + k/2 \Bigr)
\\
&&\qquad\leq\sum_{k=1}^\infty e^{(k+1)\kappa} \sum
_{j=0}^{\lceil t^2\rceil
}2\P_0 \Bigl(
\sup_{s\in[0,1/t]} B_s \geq k/2 \Bigr) \leq c_2
t^2\sum_{k=1}^\infty
e^{(k+1)\kappa-k^2t/8}.
\end{eqnarray*}
Writing out the Bessel density and applying the Markov property and
then Lem\-ma~\ref{bessbound}, and using that $z+y+1\leq z(y+2)$ for all
$z\geq1$, the second double sum is bounded above by
\begin{eqnarray*}
&&\sum_{k=1}^\infty e^{(k+1)\kappa}\sum
_{j=1}^{\lceil
(y+1)t^{4/3}\rceil} \int_{k(j/t+1)^{2/3}/2+y+1}^\infty
\frac
{ze^{-(z-y-1)^2 t/2j}}{(y+1)\sqrt{2\pi j/t}}\hat P_z(B_{t-j/t}\leq
y+2)\,dz
\\
&&\qquad\leq\sum_{k=1}^\infty e^{(k+1)\kappa}\sum
_{j=1}^{\lceil
(y+1)t^{4/3}\rceil} \int_{k(j/t+1)^{2/3}/2}^\infty
\frac
{ze^{-z^2t/2j}}{\sqrt{2\pi j/t}}\cdot\frac{\gamma
(y+2)^3}{(t-j/t)^{3/2}} \,dz
\\
&&\qquad\leq c_3\frac{(y+2)^3}{t^{3/2}}\sum_{k=1}^\infty
e^{(k+1)\kappa
}\sum_{j=1}^{\lceil(y+1)t^{4/3}\rceil} \int
_{k(j/t)^{1/6}/2}^\infty(j/t)^{1/2}
ze^{-z^2/2}\,dz
\\
&&\qquad\leq c_3\frac{(y+2)^3}{t^{3/2}}\sum_{k=1}^\infty
\sum_{j=1}^{\infty
} k(j/t)^{2/3}e^{(k+1)\kappa-k^2 j^{1/3}/8t^{1/3}}.
\end{eqnarray*}
Writing out the Bessel density and again using that $z+y+1\leq z(y+2)$
for all $z\geq1$, we see that the third double sum is bounded above by
\begin{eqnarray*}
&&\sum_{k=1}^\infty e^{(k+1)\kappa}\sum
_{j=\lceil(y+1)t^{4/3}\rceil
+1}^{\lceil t^2 - (y+1)t^{4/3}\rceil}\int_{k\kappa
(j/t)/2+y+1}^\infty
\frac{ze^{-(z-y-1)^2 t/2j}}{(y+1)\sqrt{2\pi
j/t}} \,dz
\\
&&\qquad\leq\sum_{k=1}^\infty e^{(k+1)\kappa}\sum
_{j=\lceil
(y+1)t^{4/3}\rceil+1}^{\lceil t^2 - (y+1)t^{4/3}\rceil}\int_{k\kappa
(j/t)t^{1/2}/2j^{1/2}}^\infty(j/t)^{1/2}
ze^{-z^2/2}\,dz
\\
&&\qquad\leq\sum_{k=1}^\infty c_4 k
t^{2/3} e^{(k+1)\kappa-k^2(y+1)^{1/3}t^{1/9}/8}.
\end{eqnarray*}
Finally, the fourth double sum is essentially the time reversal of the
second double sum: applying Lemma~\ref{bessbound} and the Markov
property, and then writing out the Bessel density, we see that the
fourth double sum is bounded above by
\begin{eqnarray*}
\hspace*{-4pt}&&\sum_{k=1}^\infty e^{(k+1)\kappa} \sum
_{j=\lceil t^2 -
(y+1)t^{4/3}\rceil}^{\lceil t^2\rceil-1} \int
_{y+1+(t-j/t+1)^{2/3}k/2}^\infty
\frac{\gamma z^2}{(j/t)^{3/2}} \hat P_z(B_{t-j/t}\leq y+2)\,dz
\\
\hspace*{-4pt}&&\qquad\leq\sum_{k=1}^\infty e^{(k+1)\kappa}\\
\hspace*{-4pt}&&\qquad\quad{}\times\hspace*{-1pt}
\sum_{j=\lceil t^2 -
(y+1)t^{4/3}\rceil}^{\lceil t^2\rceil-1} \int
_{y+1+(t-j/t+1)^{2/3}k/2}^\infty
\int_0^{y+2} \hspace*{-1pt}\frac{\gamma z w
e^{-(w-z)^2/2(t-j/t)}}{(j/t)^{3/2}\sqrt{2\pi(t-j/t)}}\,dw \,dz
\\
\hspace*{-4pt}&&\qquad\leq\sum_{k=1}^\infty e^{(k+1)\kappa}
\sum_{j=\lceil t^2 -
(y+1)t^{4/3}\rceil}^{\lceil t^2\rceil-1} c_5
\frac{(y+2)^2}{t^{3/2}} \\
\hspace*{-4pt}&&\qquad\quad{}\times\int_{y+1+(t-j/t+1)^{2/3}k/2}^\infty
\frac{z}{\sqrt{t-j/t}} e^{-(z-y-1)^2/2(t-j/t)} \,dz
\\
\hspace*{-4pt}&&\qquad\leq\sum_{k=1}^\infty e^{(k+1)\kappa}
\sum_{j=\lceil t^2 -
(y+1)t^{4/3}\rceil}^{\lceil t^2\rceil-1} c_6
\frac{(y+2)^3}{t^{3/2}} \int_{(t-j/t)^{1/6}k/2}^\infty(t-j/t)^{1/2}z
e^{-z^2/2} \,dz
\\
\hspace*{-4pt}&&\qquad\leq c_7\frac{(y+2)^3}{t^{3/2}}\sum_{k=1}^\infty
\sum_{j=\lceil t^2
- (y+1)t^{4/3}\rceil}^{\lceil t^2\rceil-1} k(t-j/t)^{2/3}
e^{(k+1)\kappa-k^2(t-j/t)^{1/3}/8}.
\end{eqnarray*}
For $t\geq1$ each of these terms is smaller than a constant times
$(y+2)^3 t^{-3/2}$, as required.
\end{pf}

%
\begin{prop}\label{bramsonupper}
There exists a constant $C_5$ such that
\[
\P\bigl(\exists u\in N(t)\dvtx X_u(t) \geq\sqrt{2}t -
\tfrac{3}{2\sqrt
2}\log t + y \bigr) \leq C_5(y+2)^4
e^{-\sqrt{2} y},
\]
whenever $t\geq1$ and $y\in[0,\sqrt t]$.
\end{prop}

\begin{pf}
We need to check that with high probability no particles ever go above
$\beta s + L(s) + y$ for $s\in[0,t]$. To this end define
\[
\tau= \inf\bigl\{s\in[0,t]\dvtx\exists u\in N(s) \mbox{ with } X_u(s)
> \beta s + L(s) + y\bigr\}.
\]
We claim that
\[
\E[\Gamma| \tau< t] \geq c_1
\]
for some constant $c_1>0$ not depending on $t$ or $y$; essentially if a
particle has already reached $\beta s + L(s) + y$, then it has done the
hard work, and the usual cost $e^{-\sqrt2 y}$ of reaching $\beta t$ disappears.
Choose $s<t$. On the event $\tau=s$, let $v$ be the particle at
position $\beta s + L(s) + y$ at time $s$ and define $N_v(r)$ to be the
set of descendants of particle $v$ at time $r$, for $r\geq s$. Then on
the event $\tau= s$
\begin{eqnarray*}
\Gamma&\geq&\#\bigl\{u\in N_v(t)\dvtx X_u(r) -
X_u(s) \leq\beta_s(r-s) + 1 \ \forall r\in[s,t],
\\
&&\hspace*{57.5pt}\beta_s(t-s)-1 \leq X_u(t)-X_u(s) \leq
\beta_s(t-s)\bigr\},
\end{eqnarray*}
where $\beta_s = (\beta- \frac{L(s) + y}{t-s} )\wedge
0$. It is easy to check that $\beta_s \leq\sqrt2 - \frac{3}{2\sqrt
2}\frac{\log(t-s)}{t-s}$.
Thus $\E[\Gamma|\tau=s] \geq\E[H(1,t-s)]$, and by Lemma \ref
{firstmoment}, if $s\leq t-1$, then
\[
\E[\Gamma|\tau=s] \geq c_2.
\]
If $s>t-1$, then $\E[\Gamma|\tau=s]$ is at least the probability
that a single Brownian motion $B_r, r\geq0$, remains within $[-1,1]$
for all $r\in[0,1]$, and satisfies $B_1\in[-1,0]$. This establishes
our claim, so
\[
\E[\Gamma| \tau<t] \geq c_1 \quad\mbox{and}\quad \E[\Gamma] \leq
C_4(y+2)^4 e^{-\sqrt2 y}.
\]
But then
\[
\P(\tau<t) \leq\frac{\E[\Gamma]\P(\tau<t)}{\E[\Gamma\ind_{\{
\tau<t\}}]} = \frac{\E[\Gamma]}{\E[\Gamma|\tau<t]} \leq\frac
{C_4}{c_1}(y+2)^4
e^{-\sqrt2 y}.
\]
Applying Markov's inequality, we have
\begin{eqnarray*}
\P\bigl(\exists u\in N(t)\dvtx X_u(t) \geq\sqrt{2}t -
\tfrac{3}{2\sqrt
2}\log t + y \bigr) &\leq&\P(\Gamma\geq1) + \P(\tau<t)\\
&\leq&
c_3(y+2)^4 e^{-\sqrt{2} y}
\end{eqnarray*}
as required.
\end{pf}

\begin{pf*}{Proof of Theorem~\ref{bramsonthm}}
We have shown that for $t\geq1$ and $y\in[0,\sqrt t]$, for some
constants $C_3,C_5\in(0,\infty)$,
%
%
\begin{equation}
\label{bounds} C_3 e^{-\sqrt{2}y}\leq\P\bigl(M_t >
\sqrt{2}t-\tfrac{3}{2\sqrt
2}\log t + y\bigr) \leq C_5(y +
2)^4 e^{-\sqrt{2}y}.\vadjust{\goodbreak}
\end{equation}
Thus there exists $\delta>0$ such that if we define $\tilde m(t):=
\sup\{x\in\R\dvtx\P(M_t\leq x) \leq1-\delta\}$, then
\[
\tilde m(t) = \sqrt2 t - \tfrac{3}{2\sqrt2}\log t + O(1).
\]

Fix $\varepsilon>0$. Choose $L$ such that $\E[(1-\delta)^{-|N(L)|}]
< \varepsilon/2$, and then $a$ such that $\P(M^-_L < -a) <
\varepsilon/2$ where $M^-_t = \min_{u\in N(t)} X_u(t)$ is the minimum
at time $t$. For a particle $u\in N(L)$ and $t>L$, we let $M^{(u)}_t =
\max_{v\in N(t)\dvtx u\leq v} X_v(t)$ be the maximum position among
descendants of $u$ at time $t$. Then for $t>L$,
\begin{eqnarray*}
&&
\P\bigl(M_t < \tilde m(t-L) - a\bigr)\\
&&\qquad\leq\P\bigl(M^-_L
< -a\bigr) + \P\Bigl(M^-_L \geq-a, \max_{u\in N(L)}
M^{(u)}_t < \tilde m(t-L)-a\Bigr)
\\
&&\qquad\leq\P\bigl(M^-_L < -a\bigr) + \E\bigl[\P\bigl(M_{t-L}<
\tilde m(t-L)\bigr)^{|N(L)|}\bigr]
\\
&&\qquad\leq\varepsilon/2 + \varepsilon/2 = \varepsilon.
\end{eqnarray*}
Thus $M_t - \tilde m(t)$ is tight, and we deduce that also
\[
m(t) = \sqrt2 t - \tfrac{3}{2\sqrt2}\log t + O(1).
\]
\upqed
\end{pf*}

\begin{rmk*}
It may be helpful to note that Bessel processes are not a
\textit{necessary} ingredient in our proof. One may instead replace every
appearance of a Bessel change of measure with a calculation of the
probability for a Brownian motion to stay positive, using the
reflection principle. Indeed the Bessel density can be derived directly
in this way, giving an indication that the two approaches are
interchangeable. Using the Bessel change of measure, however, conforms
with a method that works with a variety of similar problems. The
general principle is that if one wishes to calculate the number of
particles in a certain set, then one finds the martingale that forces
one particle (the \textit{spine}) to stay within that set, and studies
the corresponding measure change.
\end{rmk*}

\section{\texorpdfstring{Proof of Theorem \protect\ref{hushithm}}{Proof of Theorem 2}}

For Theorem~\ref{hushithm} we proceed via a series of four results,
each proving one of the upper or lower bounds in one of the statements
(\ref{convliminf}) or (\ref{convlimsup}).

%
\begin{lem}\label{nooneio}
The upper bound in (\ref{convliminf}) holds
\[
\liminf_{t\to\infty} \frac{M_t - \sqrt{2}t}{\log t} \leq-\frac
{3}{2\sqrt2} \qquad\mbox{almost
surely}.
\]
\end{lem}

\begin{pf}
To rephrase the statement of the lemma, we show that for any
\mbox{$\varepsilon>0$}, there are arbitrarily large times such that there are
no particles above $\sqrt{2}t - (3/2\sqrt{2} - \varepsilon)\log t$.
Choose $R>2/\varepsilon$, let $t_1 = 1$ and for $n>1$ let $t_n =
e^{Rt_{n-1}}$. Define
\[
E_n = \bigl\{\exists u\in N(t_n)\dvtx
X_u(t_n) > \sqrt{2}t_n - \bigl(
\tfrac
{3}{2\sqrt2}-\varepsilon\bigr)\log t_n\bigr\}
\]
and
\[
F_n = \bigl\{\bigl|N(t_n)\bigr|\leq e^{2t_n},
\bigl|X_u(t_n)\bigr| \leq\sqrt{2}t_n \ \forall u
\in N(t_n)\bigr\}.
\]
We know that $F_n$ happens for all large $n$, so it suffices to show that
\[
\P\biggl(\bigcap_{k\geq n} (E_k\cap
F_k) \biggr) = \lim_{N\to\infty
} \prod
_{k=n}^N \P\Biggl(E_k\cap
F_k \Big| \bigcap_{j=n}^{k-1}(E_j
\cap F_j) \Biggr) \to0 \qquad\mbox{as } n\to\infty.
\]
For a particle $u$, let $E^u_n$ be the event that some descendant of
$u$ at time $t_n$ has position larger than $\sqrt{2}t_n - (\frac
{3}{2\sqrt2}-\varepsilon)\log t_n$. Also let $s_n = t_n-t_{n-1}$ and
\begin{eqnarray*}
G_n&=&\biggl\{\exists u \in N(s_n)\dvtx\\[-4pt]
&&\hspace*{6pt}X_u(s_n) > \sqrt{2}s_n -
\frac{3}{2\sqrt2}\log s_n
+ \frac
{3}{2\sqrt2}\log\biggl(
\frac{t_n-t_{n-1}}{t_n}\biggr) + \varepsilon\log t_n\biggr\}.
\end{eqnarray*}
Then
\begin{eqnarray*}
\P\Biggl(E_k\cap F_k \Big| \bigcap
_{j=n}^{k-1}(E_j\cap F_j)
\Biggr) &\leq&
\P\Biggl(E_k \Big| \bigcap_{j=n}^{k-1}(E_j
\cap F_j) \Biggr)
\\[-2pt]
&\leq&
\P\Biggl(\bigcup_{u\in N(t_{k-1})} E^u_k
\Big| \bigcap_{j=n}^{k-1}(E_j\cap
F_j) \Biggr)
\\
&\leq&
e^{2t_{k-1}} \P(G_k)
\\
&\leq&
C_5(\log t_k+2)^4
t_k^{2/R} \biggl(1-\frac{t_{k-1}}{t_k}
\biggr)^{-3/2} t_k^{-\varepsilon},
\end{eqnarray*}
where the last inequality used Proposition~\ref{bramsonupper}. Since
we chose $R>2/\varepsilon$, this is much smaller than 1 when $k$ is
large, as required.
\end{pf}

%
\begin{lem}\label{nooneev}
The upper bound in (\ref{convlimsup}) holds
\[
\limsup_{t\to\infty} \frac{M_t - \sqrt{2}t}{\log t} \leq-\frac
{1}{2\sqrt2} \qquad\mbox{almost
surely}.
\]
\end{lem}

\begin{pf}
We show that for large $t$ and any $\varepsilon>0$, there are no
particles above $\sqrt{2}t - (1/2\sqrt2 - 2\varepsilon)\log t$. By
Proposition~\ref{bramsonupper},
\[
\P\bigl(\exists u\in N(t)\dvtx X_u(t)>\sqrt{2}t -
\bigl(\tfrac{1}{2\sqrt
2}-\varepsilon\bigr)\log t\bigr) \leq C_5(\log
t+2)^4 t^{-1-\varepsilon\sqrt2}.
\]
Thus for any lattice times $t_n\to\infty$, by Borel--Cantelli,
\[
\P\bigl(\exists u\in N(t_n)\dvtx X_u(t_n)>
\sqrt{2}t_n - \bigl(\tfrac
{1}{2\sqrt2}-\varepsilon\bigr)
\log t_n \mbox{ for infinitely many } n\bigr) = 0.\vadjust{\goodbreak}
\]
It is now a simple exercise using the exponential tightness of Brownian
motion and the fact that we may choose $t_n-t_{n-1}$ arbitrarily small
to ensure that no particle goes above $\sqrt{2}t - (\frac{1}{2\sqrt
2} - 2\varepsilon)\log t$ for any time $t$.
\end{pf}

%
\begin{lem}\label{existev}
The lower bound in (\ref{convliminf}) holds:
\[
\liminf_{t\to\infty} \frac{M_t - \sqrt{2}t}{\log t} \geq-\frac
{3}{2\sqrt2} \qquad\mbox{almost
surely}.
\]
\end{lem}

\begin{pf}
We show that for large $t$ and any $\varepsilon>0$, there are always
particles above $\sqrt{2}t - (\frac{3}{2\sqrt2} + 3\varepsilon)\log t$.
Let
\[
A_t = \bigl\{\not\exists u\in N(t)\dvtx X_u(t) > \sqrt2 t
- \bigl(\tfrac
{3}{2\sqrt2} + 2\sqrt2\varepsilon\bigr)\log t\bigr\}
\]
and
\[
B_t = \bigl\{\bigl|N(\varepsilon\log t)\bigr|\geq t^{\varepsilon/2},
X_v(\log t) \geq- \sqrt2\varepsilon\log t \ \forall v\in N(
\varepsilon\log t)\bigr\}.
\]
Define $N(v;t)$ to be the set of descendants of particle $v$ that are
alive at time $t$. Let $l_t = t - \varepsilon\log t$. Then for all
large $t$,
\begin{eqnarray*}
\P(A_t\cap B_t) &\leq& \E\biggl[\prod
_{v\in N(\varepsilon\log t)} \P\biggl(\not\exists u\in N(v;t)\dvtx\\[-5pt]
&&\qquad\quad\hspace*{39.6pt} X_u(t)>
\sqrt{2}t - \biggl(\frac
{3}{2\sqrt2}+2\sqrt2\varepsilon\biggr)\log t \Big|
\F_{\log t}\biggr) \ind_{B_t} \biggr]
\\
&\leq& \E\biggl[\prod_{v\in N(\log t)}\P\biggl(\not\exists u\dvtx
X_u(l_t) > \sqrt2 l_t -
\frac{3}{2\sqrt2}\log l_t\biggr)\ind_{B_t} \biggr]
\\
&\leq& (1-C_3)^{t^{\varepsilon/2}},
\end{eqnarray*}
where $C_3>0$ is the constant from Proposition~\ref{bramsonlower}. Thus
by Borel--Cantelli, for any lattice times $t_n\to\infty$, $\P
(A_{t_n}\cap B_{t_n} \mbox{ infinitely often})=0$. But for all large
$t$, $|N(\varepsilon\log t)|\geq e^{({\varepsilon}\log t)/{2}} =
t^{\varepsilon/2}$ and $X_v(\varepsilon\log t) \geq- \sqrt2
\varepsilon\log t$ for all $v\in N(\log t)$, so we deduce that $\P
(A_{t_n} \mbox{ infinitely often}) = 0$. Then it is again a simple task
using the exponential tightness of Brownian motion to
check that no particles move further than $(3-2\sqrt2)\varepsilon\log
t$ between lattice times infinitely often (provided that we choose $t_n
- t_{n-1}$ small enough).
\end{pf}

%
\begin{prop}\label{existio}
The lower bound in (\ref{convlimsup}) holds:
\[
\limsup_{t\to\infty} \frac{M_t - \sqrt{2}t}{\log t} \geq-\frac
{1}{2\sqrt2} \qquad\mbox{almost
surely}.
\]
\end{prop}

\begin{pf}
This is related to the proof of the lower bound in Theorem
\ref{bramsonthm}; the basic idea is similar to that in the proof given by\vadjust{\goodbreak}
Hu and Shi~\cite{hushiminimalposcritmgconvBRW}. We let $\beta_t =
\sqrt2 - \frac{1}{2\sqrt2}\frac{\log t}{t}$ and
\[
V(t) = \bigl\{v\in N(t)\dvtx X_v(r) < \beta_t r + 1 \
\forall r\leq t, \beta_t t -1 \leq X_v(t) \leq
\beta_t t \bigr\}
\]
and define
\[
I_n = \int_n^{2n}
\ind_{\{V(t)\neq\varnothing\}} \,dt.
\]
We estimate the first two moments of $I_n$. From our earlier lower
bound on $\P(H(y,t)\neq0)$ (from the proof of Proposition \ref
{bramsonlower}, taking $y = \frac{1}{\sqrt2}\log t$) we get
\[
\E[I_n] = \int_n^{2n} \P
\bigl(V(t)\neq\varnothing\bigr) \,dt \geq c_1\int_n^{2n}
e^{-(\sqrt2\cdot\log t)/{\sqrt2}} \,dt = c_2.
\]
Now,
\begin{eqnarray*}
\E\bigl[I_n^2\bigr] &=& \E\biggl[\int
_n^{2n}\int_n^{2n}
\ind_{\{V(s)\neq
\varnothing\}}\ind_{\{V(t)\neq\varnothing\}} \,ds \,dt \biggr] \\
&=& 2\int
_n^{2n}\int_n^t
\P\bigl(V(s)\neq\varnothing, V(t)\neq\varnothing\bigr) \,ds \,dt.
\end{eqnarray*}
But whenever $s\leq t$,
%
%
\begin{equation}
\label{vv} \P\bigl(V(s)\neq\varnothing, V(t)\neq\varnothing\bigr) \leq\E
\bigl[\bigl|V(s)\bigr|\bigl|V(t)\bigr| \bigr] = \E\bigl[\bigl|V(s)\bigr|
\E\bigl[\bigl|V(t)\bigr| |\F_s
\bigr] \bigr]
\end{equation}
and letting $N(u;t)$ be the set of descendants of particle $u$ that are
alive at time $t$,
\[
\E\bigl[\bigl|V(t)\bigr| |\F_s \bigr] = \sum_{u\in N(s)}
\E\biggl[ \sum_{v\in N(u;t)}\ind_{\{v\in V(t)\}}\Big|
\F_s \biggr].
\]
Now for any $s,t > 0$, let
\[
A_t(s) = \bigl\{u\in N(s)\dvtx X_u(r) < \beta_t
r + 1 \ \forall r\leq s\bigr\}
\]
and
\[
B_t(s) = \bigl\{u\in
N(s)\dvtx\beta_t s - 1 \leq X_u(s)\leq\beta_t
s\bigr\}.
\]
By the Markov property, and then applying the many-to-one lemma with
$f(r)=\beta_t(r-s)$ and $\alpha= \beta_t s - x + 1$, we have
\begin{eqnarray*}
&&\E\biggl[\sum_{v\in N(u;t)}\ind_{\{v\in V(t)\}}\Big|
\F_s \biggr]
\\
&&\qquad=\ind_{\{u\in A_t(s)\}}\\
&&\qquad\quad{}\times\E\biggl[\sum_{v\in
N(t-s)}
\ind_{\{X_v(r-s)+x<\beta_t r + 1\ \forall r\leq
t-s, \beta_t t - 1 \leq X_v(t-s)+x\leq\beta_t t\}
} \biggr] \bigg|_{x=X_u(s)}
\\
&&\qquad=\ind_{\{u\in A_t(s)\}}e^{t-s}\Q\biggl[\frac{(\beta_t s - x
+ 1)\ind_{\{\beta_t t - x - 1 \leq\xi_{t-s} \leq\beta_t t - x\}
}}{(\beta_t t - x + 1 - \xi_{t-s})e^{\beta_t \xi_{t-s} - \beta
_t^2(t-s)/2}}
\biggr]\bigg|_{x=X_u(s)}
\\
&&\qquad\leq\ind_{\{u\in A_t(s)\}} e^{t-s} \frac{\beta_t s - X_u(s) +
1}{e^{\beta_t^2 t-\beta_t X_u(s) - \beta_t - \beta_t^2(t-s)/2}}
\Q[\ind_{\{\beta_t t - x - 1 \leq\xi_{t-s}\leq\beta
_t t - x\}} ]|_{x=X_u(s)}
\\
&&\qquad\leq c_3 e^{-2s}t^{1/2}e^{({s\log t})/({2t})}\\
&&\qquad\quad{}\times
\ind_{\{u\in A_t(s)\}} \frac{\beta_t s - X_u(s) + 1}{e^{-\beta_t
X_u(s)}} \Q\bigl(\xi_t\in
B_t(t)|\xi_s = x \bigr)\big|_{x=X_u(s)},
\end{eqnarray*}
where for the last equality we used the fact that Bessel processes
satisfy the Markov property. Substituting back into (\ref{vv}) and
applying the many-to-two lemma, we get
\begin{eqnarray*}
&&\P\bigl(V(s)\neq\varnothing, V(t)\neq\varnothing\bigr)
\\
&&\qquad\leq\E\biggl[\sum_{u,v\in N(s)} \ind_{\{u\in V(s)\}}c_3
e^{-2s}t^{1/2}e^{({s\log t})/({2t})}\\
&&\hspace*{44.3pt}\qquad\quad{}\times\ind_{\{v\in A_t(s)\}} \bigl(
\beta_t s - X_v(s) + 1\bigr) e^{\beta_t X_v(s)}
\\
&&\hspace*{74pt}\qquad\quad{} \times\Q\bigl(\xi_t\in B_t(t) |\xi_s = x
\bigr)\big|_{x=X_v(s)} \biggr]
\\
&&\qquad= e^{3s}\Q\biggl[\frac{\ind_{\{T>s\}}}{\zeta^1(s)}\ind_{\{\xi
^1_s\in B_s(s)\}}
c_3 e^{-2s}t^{1/2}e^{({s\log t})/({2t})}\\
&&\qquad\quad\hspace*{60.6pt}{}\times
\zeta^1(s)e^{\beta_t^2 s/2} \Q\bigl(\xi^1_t
\in B_t(t)|\xi^1_s \bigr) \biggr]
\\
&&\qquad\quad{} + e^{2s}\Q\biggl[\frac{e^T\zeta^1(T)\ind_{\{T\leq
s\}}}{\zeta^1(s)\zeta^2(s)} \ind_{\{\xi^1_s\in B_s(s)\}}
c_3e^{-2s}t^{1/2}e^{({s\log t})/({2t})}\\
&&\qquad\quad\hspace*{113.5pt}{}\times
\zeta^2(s)e^{\beta_t^2
s/2}\Q\bigl(\xi^2_t
\in B_t(t)|\xi^2_s \bigr) \biggr]
\\
&&\qquad\leq c_4 t^{1/2}\Q\bigl(\xi^1_s
\in B_s(s), \xi^1_t\in B_t(t)
\bigr)
\\
&&\qquad\quad{} + c_4 t^{1/2} \Q\biggl[ \frac{(\beta_t r - \xi^1_r
+ 1)e^{T+\beta_t \xi^1_T - \beta_t^2 T/2}}{(\beta_t s - \xi^1_s +
1)e^{\beta_t \xi^1_s - \beta_t^2 s/2}}
\ind_{\{T\leq s\}}e^s\ind_{\{
\xi^1_s\in B_s(s), \xi^2_t \in B_t(t)\}} \biggr]
\\
&&\qquad\leq c_4 t^{1/2}\Q\bigl(\xi^1_s
\in B_s(s), \xi^1_t\in B_t(t)
\bigr)
\\
&&\qquad\quad{} + c_5 t^{1/2} e^{({s\log t})/({2t})} \Q\bigl[\bigl(
\beta_t T - \xi^1_T + 1\bigr)e^{2T - ({T\log t})/({2t}) - \beta_t(\beta
_t T -
\xi^1_T + 1)}\\
&&\qquad\quad\hspace*{151pt}\hspace*{25.8pt}{}\times
\ind_{\{T\leq s\}}\ind_{\{\xi^1_s\in B_s(s),
\xi^2_t \in B_t(t)\}} \bigr].
\end{eqnarray*}
We must now estimate the last line above. The $\Q(\cdot)$ part of the
first term is the probability that a Bessel process is near the origin
at time $s$, and then again at time~$t$; so the first term is no bigger
than a constant times $t^{1/2}s^{-3/2}(t-s+1)^{-3/2}$. Then using
notation from Section~\ref{besssec}, the expectation $\Q[\cdot]$ in
the second term is
\begin{eqnarray*}
\hspace*{-4pt}&&\hat P \bigl[Y^1_\tau e^{2\tau- ({\tau\log t})/({2t})- \beta_t
Y^1_\tau}
\ind_{\{\tau\leq s\}}\\
\hspace*{-4pt}&&\quad\hspace*{0pt}{}\times\ind_{\{(\log s)/({2\sqrt2}) -
({s\log t})/{({2\sqrt2}t)} + 1 \leq Y^1_s \leq
(\log s)/({2\sqrt2}) - ({s\log t})/({{2\sqrt2}t}) + 2\}} \ind_{\{
1 \leq Y^2_t \leq2\}} \bigr].
\end{eqnarray*}
Thus by Lemma~\ref{besslem2},
\[
\P\bigl(V(s)\neq\varnothing, V(t)\neq\varnothing\bigr) \leq c_6
\bigl(t^{-2} + t^{-1}(t-s+1)^{3/2}\bigr)
\]
and hence
\[
\E\bigl[I_n^2\bigr] \leq2 c_6 \int
_n^{2n} \int_n^t
\bigl(t^{-2} + t^{-1}(t-s+1)^{3/2}\bigr) \,ds \,dt \leq
c_{7},
\]
so
\[
\P(I_n > 0) \geq\P\bigl(I_n \geq
\E[I_n]/2\bigr) \geq\frac{\E[I_n]^2}{4\E
[I_n^2]} \geq c_8>0.
\]
When $n$ is large, at time $2\delta\log n$ there are at least
$n^{\delta}$ particles, all of which have position at least $-2\sqrt2
\delta\log n$. By the above, the probability that none of these has a
descendant that goes above $\sqrt2 s - \frac{1}{2\sqrt2}\log s -
2\sqrt2 \delta\log n$ for any $s$ between $2\delta\log n + n$ and
$2\delta\log n + 2n$ is no larger than
\[
(1-c_8)^{n^\delta}.
\]
The result follows by the Borel--Cantelli lemma since $\sum_n
(1-c_8)^{n^\delta} < \infty$.
\end{pf}

\begin{pf*}{Proof of Theorem~\ref{hushithm}}
The result is given by combining Lemmas~\ref{nooneio}, \ref
{nooneev} and~\ref{existev} and Proposition~\ref{existio}.
\end{pf*}

\section*{Acknowledgments}
The author would like to thank Ming Fang for pointing out an error in
an earlier draft, and Ofer Zeitouni for the argument required to finish
the proof of Theorem~\ref{bramsonthm} without resorting to applying
existing tightness results. Two referees also provided many helpful
comments and corrections. Finally, thanks go to Simon Harris and
Andreas Kyprianou for their valuable advice.


%

\printaddresses

\end{document}